\documentclass[a4paper,12pt]{scrartcl}

\usepackage[utf8]{inputenc}
\usepackage[english]{babel}
\usepackage{amssymb}
\usepackage{amsmath}
\usepackage{latexsym}
\usepackage{amsthm}
\usepackage{eucal}
\usepackage{graphicx}
\usepackage{bbm}
\usepackage{verbatim}
\usepackage{epigraph}
\usepackage{mathtools}
\usepackage[pdflatex]{crop}
\usepackage[bookmarks]{hyperref}
\usepackage{tikz}
\usepackage{tikz-cd}

\setkomafont{disposition}{\normalfont\bfseries}

\hypersetup{
pdfauthor={Paolo Perrone},%
pdftitle={Dual Connections},%
colorlinks, linktocpage=true, pdfstartpage=1, pdfstartview=FitV,%
breaklinks=true, pdfpagemode=UseNone, pageanchor=true, pdfpagemode=UseOutlines,%
plainpages=false, bookmarksnumbered, bookmarksopen=true, bookmarksopenlevel=1,%
hypertexnames=true, pdfhighlight=/O,%
urlcolor=black, linkcolor=black, citecolor=black, 
}

\theoremstyle{definition}
\newtheorem{deph}{Definition}[section]
\theoremstyle{remark}

\theoremstyle{plain}
\newtheorem{thm}{Theorem}[section]

\newtheorem{prop}[thm]{Proposition}

\DeclareMathOperator{\Aut}{Aut}

\DeclareMathOperator{\GL}{GL}
\DeclareMathOperator{\gl}{\mathfrak{gl}}
\DeclareMathOperator{\g}{\mathfrak{g}}

\DeclareMathOperator{\R}{\mathbb{R}}

\DeclareMathOperator{\bra}{\langle}
\DeclareMathOperator{\ket}{\rangle}

\title{Dual Connections and Holonomy}
\author{Paolo Perrone}

\begin{document}

\maketitle
\addcontentsline{toc}{subsection}{Abstract}

\begin{abstract}
Dual affine connections on Riemannian manifolds have played a central role in the field of information geometry since their 
introduction in \cite{families}. 

Here I would like to extend the notion of dual connections to general vector bundles with an inner product, in the same way as a
unitary connection generalizes a metric affine connection, using Cartan decompositions of Lie algebras.
This gives a natural geometric interpretation for the Amari tensor, as a ``connection form term'' which generates 
dilations, and which is reversed for the dual connections.
\end{abstract}

\tableofcontents
\addcontentsline{toc}{subsection}{Contents}

\newpage
\section{Cartan Decomposition}\label{cartan}

Let $\gl(n,\R)$ be the Lie algebra of $\GL(n,\R)$. If $X^*$ denotes the transpose in $\gl$, 
then the map $\theta:\gl\to\gl$ given by $\theta(x):=-x^*$ has the following properties.
\begin{enumerate}
 \item It is a (linear) isomorphism of $\gl$.
 \item It respects Lie brackets: $\theta[x,y]=[\theta x, \theta y]$.
 \item It is an involution: $\theta\circ\theta=id$.
 \item If $B$ is the Killing form, then $B_\theta(x,y):=-B(x,\theta y)$ is a positive definite 
 symmetric bilinear form.
\end{enumerate}

\begin{deph}
 A homomorphic involution of a Lie algebra $\g$ which satisfies the pro\-perties above
 is called a \emph{Cartan involution}.
\end{deph}

Any real semisimple Lie algebra admits a Cartan involution, which is unique up to inner isomorphisms.
A Cartan involution $\theta$ divides $\g$ into two eigenspaces:
\begin{equation}
 g = \mathfrak{k}\oplus\mathfrak{p}\;,
\end{equation}
where the eigenvalues are respectively $+1$ and $-1$. For $\gl$, this corresponds intuitively to the 
``direction of rotations'', and ``direction of dilations''. 
We will see that dualizing a connection reverses precisely only the second subspace.

We have that:
\begin{enumerate}
 \item $[\mathfrak{k},\mathfrak{k}]=\mathfrak{k}$,
 \item $[\mathfrak{k},\mathfrak{p}]=\mathfrak{p}$,
 \item $[\mathfrak{p},\mathfrak{p}]=\mathfrak{k}$,
\end{enumerate}
which imply that:
\begin{itemize}
 \item $\mathfrak{k}$ is a Lie subalgebra of $\g$. 
 \item $\mathfrak{p}$ in general is not, and its Lie subalgebras are all commutative (and 1-dimensional).
 \item $\mathfrak{k}$ and $\mathfrak{k}$ are orthogonal for the Killing form and for $B_\theta$.
\end{itemize}

For $\gl$, the subalgebra $\mathfrak{k}$ is precisely $\mathfrak{so}$, on which indeed $\theta$ is the identity.

In other words, we can decompose any element $x$ of $\g$ into:
\begin{equation}\label{ddec}
 x = x^+ + x^-,
\end{equation}
where:
\begin{equation}
 x^+:= \dfrac{x+\theta x}{2}, \qquad  x^-:= \dfrac{x-\theta x}{2}.
\end{equation}

There is a corresponding involution for Lie groups, which for $\GL$ corresponds to $M\mapsto (M^*)^{-1}$.
This map is called $\Theta:G\to G$ and it has the property that for $x\in\g$:
\begin{equation}
 \Theta e^x = e^{\theta x}\;.
\end{equation}

We call $K$ the subgroup of $G$ which is fixed by $\Theta$. For $\GL$ it is precisely $SO$. 
In general it is generated by exponentials of $\mathfrak{k}$.
Every element of $G$ can be expressed as a product:
\begin{equation}
 M=e^k e^p\;,
\end{equation}
where $k\in\mathfrak{k}$ and $p\in\mathfrak{p}$. Note that in general this is \emph{not} equal to
$e^{k+p}$, as $k$ and $p$ may not commute. What is true, though, is the following.

\begin{prop}\label{duexp}
 Let $e^{k+p} = e^{k'}e^{p'}$, with $k,k'\in\mathfrak{k}$ and $p,p'\in\mathfrak{p}$ possibly different.
 Then :
 \begin{equation}
  e^{k-p} = e^{k'}e^{-p'}\;.
 \end{equation}
\end{prop}

To see this, it is sufficient to notice that the quantity above is precisely:
\begin{equation}
 \Theta e^{k+p} = e^{\theta(k+p)}\;.
\end{equation}

For $\GL$, this implies that any element $M$ can be written as $M=OP$, where $O$ is orthogonal, and $P$ 
is positive definite (and given by $P=M^*M$).

For all the details on Cartan involutions and decompositions, see \cite{knapp}, Chapter VI. 

\section{Holonomy}\label{holon}

Let $V$ be a vector bundle over $X$ of rank $n$. Let $p\in X$ and let $L_p$ be the set of loops pointed at $p\in X$
equipped with the usual composition of loops (see \cite{hol}).

We can view a connection on $V$ as a smooth mapping 
$\nabla:L_p\to \Aut(V_p)$, where $V_p$ is the fiber at $p$, such that $\nabla(l)$ is the transformation
that a vector at $p$ undergoes after parallel transport along $l$. It has the following properties:
\begin{enumerate}
 \item $\nabla$ maps the trivial loop to the identity.
 \item $\nabla$ preserves composition: $\nabla(ll')=\nabla(l')\circ \nabla(l)$.
 \item If $-l$ it the inverse loop of $l$, then $\nabla(-l)=(\nabla(l))^{-1}$.
\end{enumerate}
These properties, which remind of a group homomorphism, can indeed define a homomorphism provided that a suitable 
group structure is defined on the space of loops, through quotienting. But we will not need it here. 
Even without defining a group structure for loops, the properties above imply that:

\begin{prop}
 The image of $\nabla$ is a Lie subgroup of $\Aut(V_p)$.
\end{prop}

We call such image the \emph{holonomy group} of $\nabla$ at $p$, and we denote it by $Hol_p(\nabla)$. 
If $X$ is path connected, all holonomy groups are isomorphic. In that case we drop the reference to the base point, 
and simply write $Hol(\nabla)$.
Local coordinates around $p$ give an isomorphism between $\Aut(V_p)$ and $GL(n,\R)$, so that $Hol(\nabla)$ 
is isomorphic to a subgroup of $GL$. For example:
\begin{itemize}
 \item A trivial connection has trivial holonomy group.
 \item A general affine connection may have the whole $GL$ as holonomy group.
 \item A metric connection has holonomy group isomorphic to (a subgroup of) $O(n)$.
 Different metrics (or different coordinates) yield different isomorphisms. 
 \item A connection on an oriented bundle has holonomy group isomorphic to (a subgroup of) $SL(n)$.
 A metric connection here will yield $SO(n)$.
 \item Special subgroups of $O(n)$, like for example $SU(n/2)$ for $n$ even, are the holonomy groups of the 
 so-called manifolds of special holonomies.
\end{itemize}

Let $hol(\nabla)$ be the Lie algebra of the holonomy group. Then we can express the connection $\nabla$ locally,
in suitable coordinates, as a 1-form $\omega$ with values in $hol(\nabla)$, i.e., an element of $T^*X\otimes hol(\nabla)$,
mapping linearly a tangent vector $v$ to an element $x$ of the Lie algebra. If we also choose coordinates on the fiber,
we have a mapping from tangent vectors to a subalgebra of $\gl(n)$.

For details about holonomy, the reader is referred to \cite{hol}.

\section{Dual Connection}\label{dualcon}

We can put together the results of the previous two sections, and define dual connections.
 Let $\nabla:L_p\to \Aut(V_p)$ be a connection on $V$. 
Since $\Aut(V_p)$ is isomorphic to $\gl$, it admits Cartan involutions. Let $\Theta$ be such a Cartan involution.
Then $\nabla^*:L_p\to \Aut(V_p)$ defined by $\Theta\circ\nabla$ is called the \emph{dual connection} 
with respect to the Cartan involution $\Theta$. 

Different Cartan involutions will yield different dual connections, and this is equivalent to choosing
a different inner product on $V$ (if $V$ is the tangent bundle, a Riemannian metric). This is done 
by indentifying the metric adjoint $M\mapsto M^\dagger$ with $M\mapsto \Theta(M)^{-1}$.

Let now be $\bra\;,\,\ket$ an inner product on $V$. Let $l$ be a loop at $p$. Then we can decompose the parallel
transport according to Cartan as:
\begin{equation}
 \nabla(l) = e^ke^p\;.
\end{equation}
The dual connection will instead yield (see Proposition \ref{duexp}):
\begin{equation}
 \nabla^*(l) = e^ke^{-p}\;,
\end{equation}
so that if $v,w$ are vectors of $V_p$:
\begin{align}
 \bra \nabla^*(l) v, \nabla(l) w \ket &= \bra e^ke^{-p} v, e^ke^p w \ket \\
  &= \bra  v, (e^ke^{-p})^\dagger e^ke^p w \ket \\
  &= \bra  v, \Theta(e^ke^{-p})^{-1} e^ke^p w \ket \\
  &= \bra  v, \Theta(e^{-p})^{-1}\Theta(e^k)^{-1} e^ke^p w \ket \\
  &= \bra  v, \Theta(e^{p}) e^p w \ket \\
  &= \bra  v, e^{-p} e^p w \ket = \bra  v, w \ket\;,
\end{align}
which is the property traditionally defining dual connections (see \cite{amari}).

At the Lie algebra level, the connection form $\omega^*$ of the dual connection $\nabla^*$ is obtained by $\omega$
as the mapping $\omega^*:TX\to hol(\nabla)$ given by $\theta\circ\omega$. 

Applying the decomposition \eqref{ddec} to $\omega$, we get that:
\begin{align}
 \omega &= \omega^+ + \omega^-\;, \\
 \omega^* &= \omega^+ - \omega^-\;.
\end{align}
Equivalently:
\begin{equation}
\omega^* = \omega - 2\omega^-\;. 
\end{equation}

If $\omega^-=0$ we have a metric connection (see the Ambrose-Singer holonomy theorem in \cite{hol}). 
In general $\omega^-$ measure how much our connection tends to change the 
length of the vectors, and the dual connection does the opposite. 

On Riemannian manifolds, if $\nabla$ is torsion free, then $\omega^-=0$ gives precisely the Levi-Civita connection.
In general $\omega^-$ is precisely proportional to the Amari tensor (see \cite{amari}). 
This suggests that we can generalize the concept of $\alpha$-connections to general vector bundles, by taking:
\begin{align}
 \omega^{\alpha} &= \omega^+ + \alpha\omega^-\;, \\
 \omega^{-\alpha} &= \omega^+ - \alpha\omega^-\;.
\end{align}

Moreover, this way we have a very natural geometric interpretation of the Amari tensor: it can be written as 
a 1-form with values in $\mathfrak{p}$, i.e. a subspace of $hol(\nabla)$ orthogonal to $\mathfrak{k}\equiv \mathfrak{so}$.
Intuitively, it is the part of the connection which generates dilations. 
Since every subalgebra of $\mathfrak{p}$ is 1-dimensional, the parameter $\alpha$ spans it completely. 

\section{Acknowledgements}

This idea came after the very interesting conversations with prof. Jun Zhang the University of Chicago,
whom I would like to thank, during a conference in Edinburgh.

\bibliographystyle{unsrt}
\nocite{*}
\bibliography{dual}
\addcontentsline{toc}{section}{\bibname}

\end{document}